\documentclass[11pt,a4paper]{article}

\usepackage{color}
\usepackage{amsthm}
\usepackage{subfigure}
\usepackage{mathtools}
\usepackage{authblk}
\input{Pack.sty}

\usepackage{type1cm}        
\usepackage{makeidx}         
\usepackage{graphicx}        
\usepackage{multicol}        
\usepackage[bottom]{footmisc}

\usepackage{newtxtext}       %
\usepackage{newtxmath}       
\usepackage{bm}
\usepackage{multirow, multicol, mathtools}

\newcommand{\blue}[1]{\color{black}{#1}\color{black}}

\begin{document}
\title{Reduced order methods for parametrized non-linear and time dependent optimal flow control problems, towards applications in biomedical and environmental sciences}
\author[$\sharp$]{Maria Strazzullo}
\author[$\sharp$]{Zakia Zainib}
\author[$\sharp$]{Francesco Ballarin}
\author[$\sharp$]{Gianluigi Rozza}
\affil[$\sharp$]{mathlab, Mathematics Area, International School for Advanced Studies (SISSA), Via Bonomea 265, I-34136 Trieste, Italy}
\date{}                     
\setcounter{Maxaffil}{0}
\renewcommand\Affilfont{\itshape\small}
\maketitle

\abstract{
{
We introduce reduced order methods as an efficient strategy to solve parametrized non-linear and time dependent optimal flow control problems governed by partial differential equations. Indeed, \blue{the }optimal control problems require a huge computational effort in order to be solved, most of all in physical and/or geometrical parametrized settings. Reduced order methods are \blue{a reliable } and suitable approach, increasingly gaining popularity, to achieve rapid and accurate optimal solutions in several fields, such as in biomedical and environmental sciences. In this work, we employ \blue{a }POD-Galerkin reduction approach over a parametrized optimality system, derived from \blue{the }Karush-Kuhn-Tucker conditions. The methodology presented is tested on two boundary control problems, governed respectively by (\emph{i}) time dependent Stokes equations and (\emph{ii}) steady non-linear Navier-Stokes equations.}
}
\section{Introduction}
\label{intro}
\noindent Parametrized optimal flow control problems (OFCP({\small $\boldsymbol{\mu}$})s) constrained to parametrized partial differential equations (PDE({\small $\boldsymbol{\mu}$})s) are a very versatile mathematical model which arises in several applications, see e.g. 
\cite{makinen, negri2015reduced, de2007optimal}. These problems are computationally expensive and challenging even in \blue{a }simpler non-parametrized context. The computational cost becomes unfeasible when these problems involve time dependency \cite{agoshkov2006mathematical, Stoll} or non-linearity \cite{fursikov1998boundary,gunzburger1991analysis, de2007optimal}, in addition to physical and/or geometrical parametrized settings that describe several configurations and phenomena. A suitable strategy to lower this expensive computational effort is to employ reduced order methods (ROMs) in the context of \ocp s, which recast them in a cheap, yet reliable, low dimensional framework \cite{hesthaven2015certified,RozzaHuynhManzoni2013}.
We exploit \blue{these }techniques in order to solve boundary \ocp s on a bifurcation geometry \cite{rozza2012reduction} which can be considered as (\emph{i}) a riverbed in environmental sciences and as (\emph{ii}) a bypass graft for cardiovascular applications. 
\noindent In the first research field, reduced parametrized optimal control framework (see e.g. \cite{quarteroni2005numerical, quarteroni2007reduced}) can be of utmost importance. It perfectly fits in forecasting and data assimilated models and it could be exploited in order to prevent possibly dangerous natural situations \cite{StrazzulloBallarinMosettiRozza2018}. The presented test case is governed by time dependent Stokes equations, which are an essential tool in marine sciences in order to reliably simulate evolving natural phenomena.\\
Furthermore, discrepancies between computational modelling in cardiovascular mechanics and reality usually ought to high computational cost and lack of optimal quantification of boundary conditions, especially the outflow boundary conditions. In this work, we present application of the aforementioned numerical framework combining OFCP($\mu$) and reduced order methods in the bifurcation geometry. The aim is to quantify the outflow conditions automatically while matching known physiological data for different parameter-dependent scenarios \cite{ballarin2016fast, ZainibEtAl2019}. In this test case, Navier-Stokes equations will model the fluid flow.

 \noindent The work is outlined as follows: in section \ref{problem_and_methodology}, the problem formulation and the methodology are summarized. Section \ref{results} shows the numerical results for the two boundary \ocp s, \blue{based on \cite{negri2015reduced, rozza2012reduction}}. Conclusions follow in section \ref{conclusions}.

\section{Proper Orthogonal Decomposition for \ocp s}
\label{problem_and_methodology}
\noindent {In this section, we briefly describe the problem and the adopted solution strategy for time dependent non-linear boundary \ocp s: in the cases mentioned in section \ref{intro}, the reader shall take the non-linear term and time-dependent terms to be zero accordingly  \cite{Strazzullo2019, ZainibEtAl2019}. The goal of \ocp s is to find a minimizing solution for a quadratic cost functional {\small $\mathcal J$} under a PDE($\bmu$) constraint thanks to an external variable denoted as \emph{control}. {In the next section, we will show numerical results over a bifurcation geometry {\small $\Omega$} with physical and/or geometrical parametrization \blue{represented by } {\small $\pmb{\mu}\in \mathcal{D}\subset \mathbb{R}^d, d\in \mathbb{N}$}. Thus, considering the space-time domain 
{ \small $Q = \Omega \times [0,T]$} with a sufficiently regular spatial boundary 
{ \small $\partial \Omega$} \footnote{For the steady case, {\small $T = 0$} and
{ \small $Q \equiv \Omega$}}, let us define the Hilbert spaces} { \small $S = V\times P$},{ \small $Z = Z_V\times Z_P$} and { $U$} {for} state and adjoint velocity and pressure, and control variables denoted by { \small $\pmb{s} = \left( \pmb{v}, p\right)\in S$}, {\small $\pmb{u}\in U$} and { \small $\pmb{z} = \left( \pmb{w}, q\right)\in Z$}, respectively. The stability and uniqueness of the optimal solution will be guaranteed if 
{ \small $S\equiv Z$}, which will be our assumption in this work. We introduce { \small $X = S\times U$} such that {\small $\mathbf{x} = \left( \pmb{s}, \pmb{u}\right)\in X$}. Then, the problem reads: {\it given
{ \small $\pmb{\mu}\in \mathcal{D}$} , find {\small $\left(\mathbf{x}, \pmb{z}; \pmb{\mu}\right)\in X\times S$} s.t.:}
\begin{equation}
\label{optimality_system}
\begin{cases}
\begin{matrix}
\mathcal{A}\left(\mathbf{x}, \mathbf{y}; \pmb{\mu}\right) + \mathcal{B}\left(\pmb{z}, \mathbf{y}; \pmb{\mu}\right)\ +\\ \mathcal{E}\left(\pmb{v}, \pmb{w}, \mathbf{y}_{\pmb{v}}; \pmb{\mu}\right) + \mathcal{E}\left(\pmb{w}, \pmb{v}, \mathbf{y}_{\pmb{v}}; \pmb{\mu}\right)
\end{matrix} = \int_0^T \left\langle \mathcal{H}\left(\pmb{\mu}\right), \mathbf{y}\right\rangle dt,&\ \forall\ \mathbf{y}\in X,\\
\mathcal{B}\left(\mathbf{x}, \pmb{\kappa}; \pmb{\mu}\right) + \mathcal{E}\left(\pmb{v}, \pmb{v}, \pmb{\kappa}_{\pmb{w}}; \pmb{\mu}\right)  = \int_0^T \left\langle \mathcal{G}\left(\pmb{\mu}\right), \pmb{\kappa} \right\rangle dt,&\ \forall\ \pmb{\kappa}\in S,
\end{cases}
\end{equation}
where 
{\small $\mathbf{y} = \left(\mathbf{y}_{\pmb{v}}, \mathbf{y}_{p}, \mathbf{y}_{\pmb{u}}\right)$}
and {  \small $\pmb{\kappa} = (\pmb{\kappa}_{\pmb{w}}, \pmb{\kappa}_{q})$}.
{The bilinear forms
 { \small $\mathcal{A}: X\times X\rightarrow \mathbb{R}$} and {\small $\mathcal{H}\left(\pmb{\mu}\right)$} are related to the minimization cost functional {\small $\mathcal{J}$}, while the bilinear form 
{ \small $\Cal B : X \times Z \rightarrow \mathbb R$} represents the linear part of state-constraints and {\small $\mathcal{E}$} is the non-linear convection term, which will be zero for the time dependent linear case.}\\
In order to solve the optimality system \eqref{optimality_system} we exploit Galerkin Finite Elements (FE) {\it snapshots-based} {\it Proper orthogonal decomposition} ({POD})--Galerkin (see \cite{hesthaven2015certified}), summarized in table \ref{table1 :algorithm}, where the number of time steps are denoted by {\small $\mathcal{N}_t$}.  
\begin{table}[H]
\footnotesize
\begin{center}
\begin{tabular}{ll}
\hline
\hline
Offline phase: \hspace{5cm} Input: $\mu_1$ for lifting, $N$, $n$.&\\
\hspace{6.7cm} Output: Reduced order solution spaces.&\\
 \hline
$1.$ Compute snapshots $\bm{\delta}_{\mathcal{N}_{\bm{\delta}}\times \mathcal{N}_t}\left(\bm{\mu}^n\right)$ for $1\leq n\leq |\Lambda |,\ \bm{\delta} = \bm{v}, p, \bm{u}, \bm{w}, q$ and state and adjoint\\ supremizers. The global dimension of FE space discretization is {\small $\Cal N = 2\Cal N_{\boldsymbol v} + 
2\Cal N_{p} + \Cal N_{\boldsymbol u}$}.\\
$2.$ Solve eigenvalue problems $\mathbb{A}^{\bm{\delta}} \bm{\rho}^{\bm{\delta}}_{n} =  \lambda_{n}^{\bm{\delta}} \bm{\rho}^{\bm{\delta}}_{n},\ n = 1, \cdots, |\bm{\Lambda}|$, where $\mathbb{A}^{\bm{\delta}}$ is
correlation\\ matrix of snapshots.\\
$3.$ If relative energy of eigenvalues is greater than $1 - \epsilon_{tol},\ 0 < \epsilon_{tol} \ll 1$, keep corresponding\\ eigenvalue-eigenvector pairs $\left(\lambda^n,  {\bm{\rho}}^n\right)$.\\
$4.$ Construct orthonormal {POD} basis from the retained $N$ eigenvectors and add the {POD} modes\\ of the supremizers to state and adjoint velocities.\\
\hline \hline
Online phase: \hspace{5cm} Input: Online parameter $\bm{\mu}\in \mathcal{D}$.&\\
\hspace{6.7cm} Output: Reduced order solution.&\\
 \hline
 $1.$ Perform Galerkin projection to calculate reduced order coefficients such that { \small $\pmb{\delta} \approx {\mathfrak{X}}_{\pmb{\delta}}\pmb{\delta}_{\mathrm{N}} $}\\ where, { \small $\mathfrak{X}_{\pmb{\delta}}$} denotes reduced bases matrices containing all the time instances.\\
 $2.$ Solve the reduced order version of the optimality system \eqref{optimality_system}.\\
 \hline
 \hline
\end{tabular}
\end{center}
\caption{Algorithm: {POD}--Galerkin for \ocp}
\label{table1 :algorithm}
\end{table}


\noindent In order to guarantee the efficiency of \blue{the } POD--Galerkin approach, we rely on the affine assumption over the forms, i.e. every form can be written as a linear combination of $\small \bmu -$dependent functions and $\small \bmu -$independent quantities.
In this way, the system resolution is divided into parameter independent ({\it offline}) and dependent ({\it online}) phases (see table \ref{table1 :algorithm} for details) such that the expensive calculations are absorbed in the former stage and only \textit{online} stage is repeated every time the parameter 
{\small $\pmb \mu$} changes. 
From the perspective of the problem stability, to  ensure uniqueness of pressure at the reduced order level, we enrich the state and adjoint velocity space with {\it supremizers} and, to guarantee the fulfillment of Brezzi's inf-sup condition \cite{brezzi} at the reduced level, we use aggregated equivalent state and adjoint spaces. Thus, dimension of the reduced problem reduces from 
{ \small $\mathcal{N}\times\mathcal{N}_t$} to 
{ \small $13N$}.

\section{Results}
\label{results}
\subsection{Linear Time Dependent \ocp $\:$ governed by Stokes Equations}
 In this section, \blue{inspired by \cite{negri2015reduced, rozza2012reduction}}, we propose an OFCP({\small $\pmb \mu$}) governed by a time dependent Stokes equation.
First of all, 
let us introduce the smooth domain {\small $\Omega(\mu_2)$}. The parameter stretches the length of the reference domain shown in figure \ref{domain_bifurcation}, which will be indicated with {\small $\Omega$} from now on. We want to recover a measurement 
{ \small $ \pmb v_d(\mu_3) \in L^2(0,T; [L^2(\Omega)]^2)$} over the one dimensional observation domain { \small $\Gamma_{OBS}$} controlling the Neumann flux over {\small $\Gamma_C$}, with the inflow {\small $\bm{v}_{{in}}\left(\bm{\mu}\right) = (10\left(x_2 - 1)\left( 1 - x_2\right), 0\right)$}.
The setting is suited for environmental applications: we control the flow in order to avoid potentially dangerous situations in an hypothetical \emph{real time} monitoring plan on the domain, which can represent a riverbed.
\noindent The space-time domain is { \small $Q = \Omega \times [0,1]$}. Let us consider the following function spaces: 
{ \small $V = L^2(0,T; [H^1_{\Gamma D}(\Omega)]^2) \cap H^1(0,T; [H^1_{\Gamma D}(\Omega)\dual ]^2)$}, { \small $P = L^2(0,T; L^2(\Omega))$}
and 
{ \small $U = L^2(0,T; [L^2(\Omega)]^2)$} for state and adjoint velocity, state and adjoint pressure and for control, respectively. Then, we define { \small $X = (V \times P)\times U$}. 
For a given 
{ \small $\bmu \in \Cal D = [0.01, 1]\times[1, 2] \times  [0.01, 1]$}, we want to find the solution of time dependent Stokes equations which minimizes:
{
\begin{equation}
\small
\label{functional_maria}
\Cal J :=  \frac{1}{2} \int_{0}^T \int_{\Gamma_{OBS}}(\pmb v(\bmu) - \pmb v_d(\mu_3))^2 dsdt 
+  \frac{\alpha_1}{2}
\int_0^T \int_{\Gamma_C}\pmb u(\bmu)^2 dsdt 
+  \frac{\alpha_2}{2} \int_0^T \int_{\Gamma_C}
|\nabla \pmb u(\bmu) \pmb t|^2 dsdt, 
\end{equation}
} \noindent where {\small $\alpha_1 = 10^{-3}, \alpha_2 = 10^{-4}$} and {\small $\pmb t$} is the unit tangent vector to {\small  $ \Gamma_C$} and { \small $\pmb v_d(\mu_3) =
[\mu_3 (8( x_2^3 - x_2^2 -x_2 + 1) + 2(-x_2^3 - x_2^2 + x_2 + 1)), 0]$}. The cost functional penalizes not only the magnitude of the control, but also its rapid variations over the boundary.
\noindent The constrained minimization problem { \eqref{functional_maria}} is equivalent to the resolution of problem { \eqref{optimality_system}} where the considered forms are defined by:
\begin{equation*}
{\small
\begin{aligned}
\hspace{-3mm} \small \Cal A (\mathbf x,\mathbf y) & =  \int_{0}^T\int_{\Gamma_{OBS}}\pmb v \cdot \mathbf y_{\pmb v} \; dsdt  +\alpha_1 \int_0^T \int_{\Gamma_C}\pmb u \cdot \mathbf y_{\pmb u} \;  dsdt + \alpha_2 \int_0^T \int_{\Gamma_C}\nabla  \pmb u 
\pmb t \cdot \nabla \mathbf y_{\pmb u} 
\pmb t \;  dsdt, \\
\Cal B (\mathbf x, \boldsymbol{z}; \bmu) &= \int_Q{\dt{\pmb v}\cdot \boldsymbol w \;  dxdt} + \mu_1 \int_Q{\nabla{\pmb v} \cdot \nabla{\boldsymbol w} \;  dxdt} - \int_Q{p ( \nabla \cdot \boldsymbol w(\bmu))} \;  dxdt
\\ & \qquad\qquad \qquad -\int_Q{q (\nabla \cdot {(\boldsymbol v(\bmu))} \;  dxdt} - \displaystyle \int_0^T \int_{\Gamma_C} {\pmb u \cdot \boldsymbol w \; dsdt},\\
\la \Cal H(\bmu), \mathbf y \ra & =  \int_{\Gamma_{OBS}}\pmb v_d(\mu_3) \cdot \mathbf y_{\pmb v} \; ds, \qquad
\la \Cal G(\bmu),{q} \ra = 0,  \;\;\; \forall \mathbf y \in  X,   \vspace{-1.7cm}\\
\end{aligned}
}
\end{equation*} 

\noindent for every $ \mathbf x,\mathbf y \in X $ and $\boldsymbol \kappa \in S$. We built the reduced space with { \small $N = 35$} over a training set of {\small $70$} snapshots of global dimension {\small $131400$}, for {\small $\Cal N_t = 20$}. In time dependent applications, ROMs are of great advantage: in table \ref{speedup_bif} the
speedup index is shown with respect to $N$. The speedup represents how many ROM systems one can solve in the time of a FE
simulation. Nevertheless, we do not pay in accuracy as figure \ref{error_bif} \color{black}{and figure \ref{results_bif} show}\color{black}{}: it represents the relative error between FE and ROM variables. The relative error between FE and ROM { \small $\Cal J$} is presented in table \ref{speedup_bif}\\

\begin{table}[H]{
\begin{center}
\caption {Speedup analysis and relative error $\Cal J$.} \label{speedup_bif} 
\vspace{1mm}
\begin{tabular}{l|l|l|}
$N$  &   $ \;$Speedup   & $\;$ Relative error $\Cal J$  \\ \hline
$15$ & $\; 66338$ & $\qquad 10^{-7}$ \\
$20$ & $\; 47579$ & $\qquad 10^{-8}$ \\
$25$ & $\; 34335$ & $\qquad 10^{-8}$ \\
$30$  & $\; 22477$ & $ \qquad 10^{-9}$ \\ 
$35$  & $\; 17420$ & $ \qquad 10^{-10}$ \\ 
\end{tabular}
\end{center}
}
\end{table}


\hspace{-1cm}\begin{minipage}{.5\textwidth}
\begin{figure}[H]
\begin{center}
\includegraphics[scale = 0.3]{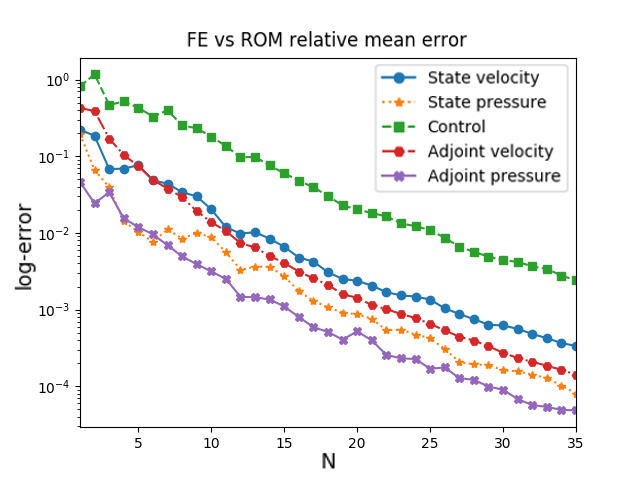}\vspace{-3mm}
\hspace{-1cm}\caption{FE vs ROM mean relative error over 50 parameters.}
\label{error_bif}
\end{center}
\end{figure}

\end{minipage}
\hspace{1cm}\begin{minipage}{.33\textwidth}
\begin{figure}[H]
\begin{center}
\hspace{7mm}\hspace{-1cm}\includegraphics[scale = 0.2]{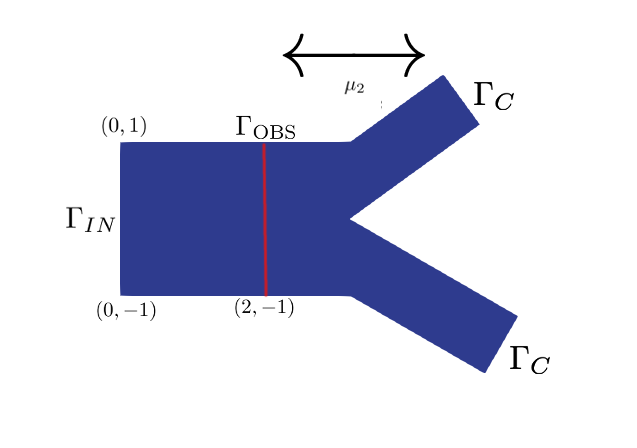}
\begin{center}
\vspace{-5mm}
\caption{Physical domain.}
\label{domain_bifurcation}
\end{center}
\end{center}
\end{figure}
\end{minipage}

\begin{figure}[H]
\begin{center}
\includegraphics[scale = 0.25]{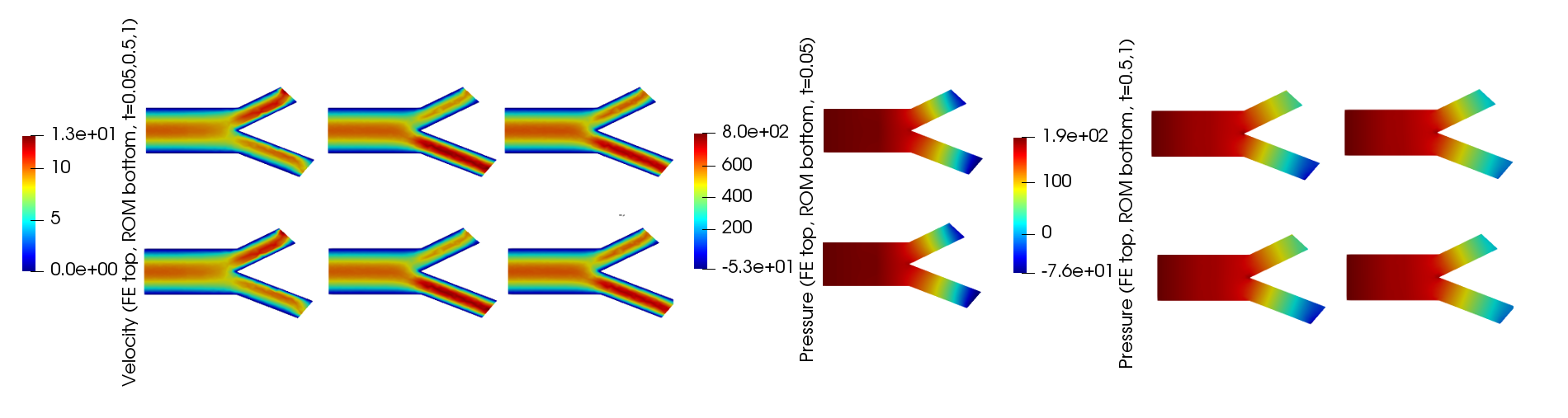}\vspace{-3mm}
\hspace{-1cm}\caption{FE (top) vs ROM (bottom) comparison of state velocity and state pressure, for t=0.05,0.5,1 and $\bmu = (0.5, 1.5, 1)$.}
\label{results_bif}
\end{center}
\end{figure}


\subsection{Non-linear steady \ocp $\:$ governed by Navier-Stokes equations}

In this section, we will demonstrate the numerical results for second test case with optimal boundary control problem governed by non-linear incompressible steady Navier-Stokes equations. We consider a bifurcation domain {\small $\Omega$} as employed in the previous example (see figure \ref{fig1: domain}), which can be considered as an idealized model of arterial bifurcation in cardiovascular problems \blue{\cite{negri2015reduced, rozza2012reduction, ZainibEtAl2019}}. Fluid shall enter the domain from {\small $\Gamma_{in}$} and shall leave through the outlets {\small $\Gamma_c$}.  In this example, physical parameterization is considered for the inflow velocity given by {\small $\bm{v}_{{in}}\left(\bm{\mu}\right) = 10{\mu}_1\left(x_2\left( 2 - x_2\right), 0\right)$} and the desired velocity, denoted by {\small $\boldsymbol{v}_d\in L^2\left(\Omega\right)$} and prescribed at the 1-D observation boundary {\small $\Gamma_{obs}$} through the following expression:
\[
\small
\bm{v}_d\left(\bm{\mu}\right) = \begin{pmatrix}
10\mu_1\left(0.8\left(\left( x_2 - 1\right)^3- \left( x_2 - 1\right)^2 -\left( x_2 - 1\right) + 1 \right) + 0.2\left(-\left( x_2 - 1\right)^3- \left( x_2 - 1\right)^2 +  x_2 \right) \right)\\
0
\end{pmatrix}.\]

The cost-functional {\small $\mathcal{J}$} is defined as:
 \begin{equation}
 \small
\label{J:NS}
\mathcal{J}\left(\bm{v}, \bm{u}; \bm{\mu}\right)=\frac{1}{2} \|\bm{v}\left(\bm{\mu}\right) - \bm{v}_d\left(\bm{\mu}\right)\|^2_{L^2\left(\Gamma_{obs}\right)} + \frac{\alpha}{2} \|\bm{u}\left(\bm{\mu}\right)\|^2_{L^2\left(\Gamma_c\right)}\\ + \frac{0.1\alpha}{2} \|\nabla\bm{u}\left(\bm{\mu}\right)\bm{t}\|^2_{L^2\left(\Gamma_c\right)}, 
\end{equation}
where {\small $\bm{t}$} is the tangential vector to {\small $\Gamma_c$}. The mathematical problem reads: {\it Given {\small $\bm{\mu}\in \mathcal{D} = \left[ 0.5, 1.5\right]$}, find {\small $\left(\bm{v}\left(\bm{\mu}\right), p\left(\bm{\mu}\right), \bm{u}\left(\bm{\mu}\right)\right)$} that minimize {\small $\mathcal{J}$} and satisfy the Navier-Stokes equations with {\small $\bm{v}_{{in}}\left(\bm{\mu}\right)$} prescribed at the inlet {\small $\Gamma_{{in}}$}, no-slip conditions at the walls {\small $\Gamma_w$} and {\small $\bm{u}\left(\bm{\mu}\right)$} implemented at {\small $\Gamma_c$} through Neumann conditions.}

At the continuous level, we consider {\small $X\left(\Omega\right)=  H^1_{\Gamma_{in}\cup\Gamma_w}\left(\Omega\right)\times L^2\left(\Omega\right)\times L^2\left(\Gamma_c\right)$}, where
\begin{equation*}
\small
H^1_{\Gamma_{in}\cup\Gamma_w}\left(\Omega\right) = \left[\bm{v}\in \left[H ^1\left(\Omega\right)\right]^2 : \bm{v} |_{\Gamma_{{in}}} = \bm{v}_{{in}}\ \text{and}\ \bm{v}|_{\Gamma_w} = 0\right].
\end{equation*}
Thus,
{\footnotesize
\begin{align*}
\mathcal{A}\left(\mathbf{x}, \bm{y}\right) & = \int_\Omega \bm{v}\left(\bm{\mu}\right)\cdot\bm{y}_{\bm{v}} d\Omega + \alpha\int_{\Gamma_c} \bm{u}\left(\bm{\mu}\right)\cdot\bm{y}_{\bm{u}} d\Gamma_c + \frac{\alpha}{10}\int_{\Gamma_c}\left( \nabla\bm{u}\left(\bm{\mu}\right)\right)\bm{t}\cdot\nabla\left(\bm{y}_{\bm{u}}\right)\bm{t} d\Gamma_c,\\
 \mathcal{B}\left(\mathbf{x}, \bm{z}\right) & = \eta\int_\Omega \nabla\bm{v}\left(\bm{\mu}\right)\cdot\nabla\bm{w}d\Omega - \int_\Omega p\left(\bm{\mu}\right)\left(\nabla\cdot \bm{w}\right) d\Omega -  \int_\Omega q\left(\nabla\cdot \bm{v}\left(\bm{\mu}\right)\right)d\Omega \\
 & - \int_{\Gamma_c} \bm{u}\left(\bm{\mu}\right)\cdot\bm{w} d\Gamma_c,\\
 \mathcal{E}\left(\bm{v}, \bm{v}, \bm{w}\right)&= \int_\Omega \left(\bm{v}\left(\bm{\mu}\right)\cdot\nabla\right)\bm{v}\left(\bm{\mu}\right)\cdot\bm{w} d\Omega\quad \text{and}\quad  \left\langle\mathcal{H}\left(\bm{\mu}\right), \mathbf{y}\right\rangle = \int_\Omega \bm{v}_d\left(\bm{\mu}\right)\cdot\bm{y}_{\bm{v}} d\Omega.
\end{align*}}
\begin{center}
\hspace{1cm}\begin{figure}[H]
\begin{minipage}{0.4\columnwidth}
\begin{center}
\includegraphics[width = 0.72\textwidth, keepaspectratio]{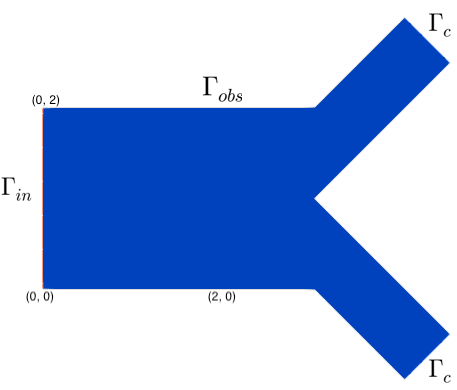}
\caption{Domain {\small $\left(\Omega\right)$}.}
\label{fig1: domain}
\end{center}
\end{minipage}\hspace{0.5cm}
\begin{minipage}{0.45\columnwidth}
\begin{center}
\begin{tabular}{|cc|c|}
\hline\hline
\multicolumn{2}{|c|}{Mesh size} & {\small $5977$}\\
\hline
\multicolumn{2}{|c|}{No. of reduced order bases {\small $N$}} & {\small $131$}\\
\hline
\multicolumn{2}{|c|}{{\small $\mathcal{D}$}} & {\small $\left[0.7, 1.5\right]$}\\
\hline
\multicolumn{2}{|c|}{{\small $| \Lambda |$}} & {\small $100$}\\
\hline
\multicolumn{2}{|c|}{offline phase} & {\small $4.9\times 10^{3}$} seconds \\\cline{2-3}
\multicolumn{2}{|c|}{online phase} & {\small $9\times 10^{1}$} seconds\\\hline
\hline

\end{tabular}
\caption{{Computational details of {POD}--Galerkin for Navier-Stokes constrained \ocp}.}
\label{tableForComputationalPerfomance}
\end{center}
\end{minipage}
\end{figure}
\end{center}

\begin{figure}
\begin{minipage}{0.44\columnwidth}
\begin{center}
\vspace{-3mm}
\includegraphics[width = 0.9\textwidth, keepaspectratio]{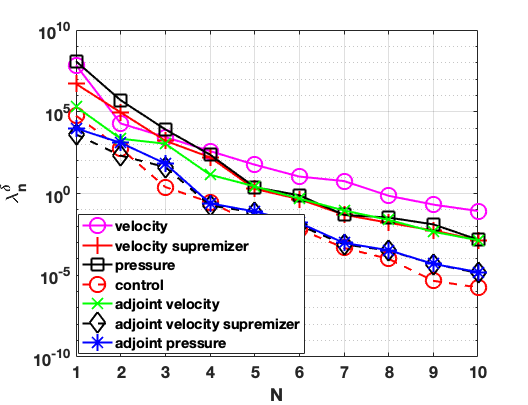}
\caption{Eigenvalues of {\small $\mathrm{N} = 10$} {POD} modes.}
\label{fig2:eigenvalues reduction}
\end{center}
\end{minipage}\hspace{0.5cm}
\begin{minipage}{0.45\columnwidth}
\begin{center}
\includegraphics[width = 0.82\textwidth, keepaspectratio]{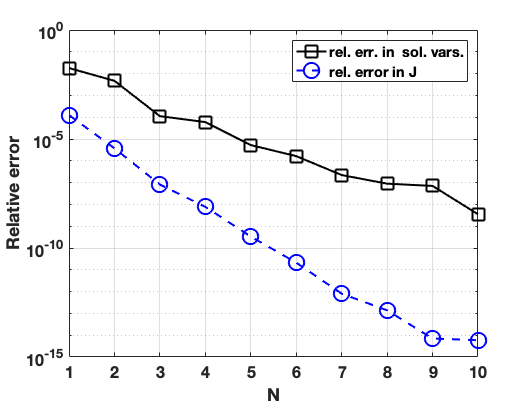}
\caption{Relative error for solution variables and {\small $\mathcal{J}$}.}
\label{fig:rel_err}
\end{center}
\end{minipage}
\end{figure}
\begin{figure}[H]
\begin{tabular}{lllll}
\multicolumn{1}{c}{\multirow{-4}{*}{{\small $|\bm{v}|$}}} & \multicolumn{1}{c}{\hspace{0.7cm}\includegraphics[width = 0.21\textwidth, keepaspectratio]{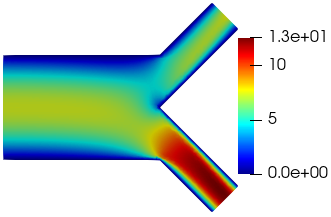}} & \multicolumn{1}{c}{\hspace{0.7cm}\includegraphics[width = 0.21\textwidth, keepaspectratio]{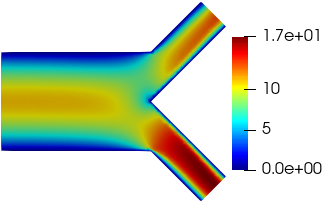}} & \multicolumn{1}{c}{\hspace{0.7cm}\includegraphics[width = 0.21\textwidth, keepaspectratio]{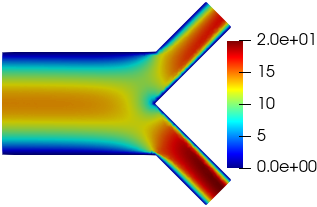}}\\
\multicolumn{1}{c}{\multirow{-4}{*}{{\small $|\bm{u}|$}}} & \multicolumn{1}{c}{\hspace{0.7cm}\includegraphics[width = 0.2\textwidth, keepaspectratio]{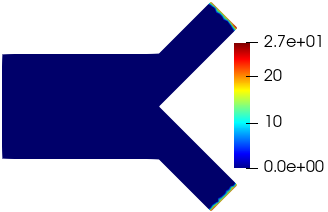}} & \multicolumn{1}{c}{\hspace{0.7cm}\includegraphics[width = 0.2\textwidth, keepaspectratio]{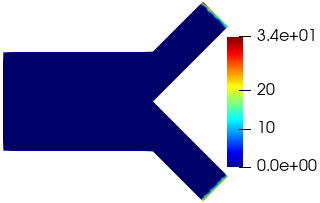}} & \multicolumn{1}{c}{\hspace{0.7cm}\includegraphics[width = 0.2\textwidth, keepaspectratio]{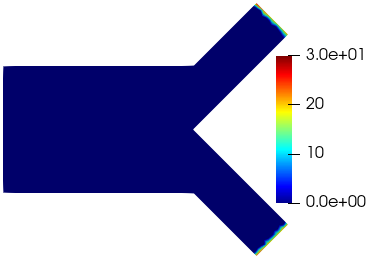}}
\end{tabular}
\caption{State velocity and control for {\small $\mu_1 = 0.7, 1.1, 1.4$}.}
\label{fig3: State velocity and control }
\end{figure}
To construct the reduced order solution spaces, we consider a sample {\small $\left(\Lambda\right)$} of {\small $100$} parameter values and solving the problem \eqref{optimality_system} through Galerkin Finite Element method, we construct the snapshot matrices for the solution variables {\small $\bm{v}, p, \bm{u}, \bm{w}, q$}. For {\small $N = 10$}, eigenvalues energy of the state, control and adjoint variables is demonstrated in figure \ref{fig2:eigenvalues reduction}. Evidently, {\small $N$} eigenvalues capture {\small $99.9\%$} of the Galerkin FE discretized solution spaces and the reduced order spaces are thus built with dimensions {\small $13\mathrm{N} +1 = 131$}. 
The state velocity and control for {\small $\mu = 0.7, 1.1, 1.4 $} are shown in figure \ref{fig3: State velocity and control }. Furthermore, we report the accumulative relative error for the solution variables and the relative error for {\small $\mathcal{J}$} in figure \ref{fig:rel_err}. The former decreases upto {\small $10^{-8}$} along with the latter decreasing upto {\small $10^{-14}$}.

\section{Concluding Remarks}
\label{conclusions}

In this work, we propose ROMs as a suitable tool to solve a parametrized boundary \ocp s for time dependent Stokes equations and steady Navier-Stokes equations. The framework proposed is suited for several \emph{many query} and \emph{real time} applications both in environmental marine sciences and bio-engineering. The reduction of the KKT system is performed through a POD-Galerkin approach, which leads to accurate surrogate solutions in a low dimensional space. This work aims at showing how ROMs can have an effective impact in the management of parametrized simulations for social life and activities, such as coastal engineering and cardiovascular problems. Indeed, the proposed framework deals with faster solving of parametrized optimal solutions which can find several applications in monitoring planning both in marine ecosystems and patient specific geometries.

\section*{Acknowledgements}
We acknowledge the support by European Union Funding for Research and Innovation -- Horizon 2020 Program -- in the framework of European Research Council Executive Agency: Consolidator Grant H2020 ERC CoG 2015 AROMA-CFD project 681447 ``Advanced Reduced Order Methods with Applications in Computational Fluid Dynamics'' (P.I. Prof. G. Rozza). We also acknowledge the INDAM-GNCS project ``Advanced intrusive and non-intrusive model order reduction techniques and applications''.
The computations in this work have been performed with RBniCS library, developed at SISSA mathLab, which is an implementation in FEniCS of several reduced order modelling techniques; we acknowledge developers and contributors to both libraries.

\bibliographystyle{plain}
\bibliography{0_OFCP_Strazzullo_Zainib}
\end{document}